\magnification=1200
\overfullrule=0mm
\baselineskip=15pt

\def\mapright#1{\smash{\mathop{\longrightarrow}\limits^{#1}}}

\vglue 0.8cm

\centerline {\bf Extensions panach\'ees autoduales.}

\bigskip

\centerline {Daniel BERTRAND  \footnote{$^{(*)}$} {Adresse de l'auteur : Institut de Math\'ematiques de Jussieu ; bertrand@math.jussieu.fr  Mots clefs : extensions panach\'ees; cat\'egories tannakiennes; repr\'esentations  unipotentes. 
 Classification AMS : 20 G 05, 20 L 05.}}

\medskip
\centerline
{{\it Novembre 2010} \footnote{$^{(1)}$}  
{On trouvera une premi\`ere version, non publi\'ee, de ce texte dans [B3]. }}

\bigskip

\bigskip
\noindent
{\bf Abstract} : {\it we study self-duality of Grothendieck's blended extensions in the context of a tannakian category. The set of equivalence classes of symmetric, resp. antisymmetric, blended extensions is naturally endowed with a torsor structure, which enables us to compute the unipotent radical of the associated  monodromy groups in various situations.}

\bigskip

La notion d'extension panach\'ee  dans une cat\'egorie
ab\'elienne 
 a \'et\'e introduite par Grothendieck [G], en liaison avec la
construction d'accouplements de monodromie relatifs 
\`a une vari\'et\'e ab\'elienne semi-stable. Nous en rappelons la
d\'efinition au
\S 1. 

\medskip

Lorsque la cat\'egorie ambiante est tannakienne, on peut parler du dual
d'une extension panach\'ee. C'est encore une extension panach\'ee, et
nous d\'ecrivons au \S 2 les types d'autodualit\'e qu'elle peut
pr\'esenter. M\^eme s'il est proche de celui 
des biextensions, notre point de vue reste lin\'eaire. Il met en valeur
la structure  de torseur 
dont est muni l'ensemble des classes d'isomorphisme d'extensions panach\'ees autoduales: voir
le Th\'eor\`eme 1 du \S 2. Les aspects bilin\'eaires sont ainsi ramen\'es
\`a des calculs
\'el\'ementaires (\S 2, Lemme 5) sur la dualit\'e dans les extensions
ordinaires.

\medskip

Une fois \'etabli le caract\`ere autodual 	d'une extension
panach\'ee $M$, elle acquiert automatiquement un signe (\S 2, Lemme
4), et ce signe permet de cerner la partie unipotente du ``groupe  de
monodromie" $G_M$ auquel $M$ donne  naissance. C'est le th\`eme du
\S 3, o\`u nous r\'eduisons la question \`a la description, classique,
de certains sous-groupes paraboliques du groupe orthogonal ou
symplectique. Pour un \'enonc\'e pr\'ecis, voir le Th\'eor\`eme 2 du \S 3, ainsi que son corollaire, qui recouvre \`a la fois les r\'esultats de Ribet [R] sur les repr\'esentations
$\ell$-adiques attach\'ees aux 1-motifs, et ceux de [B2] sur les groupes
de Galois de certaines \'equations diff\'erentielles.

\bigskip

\noindent
{\bf \S 1.  Rappels sur les extensions panach\'ees}
\bigskip

   Soient $\bf T$ une cat\'egorie ab\'elienne, et $A, B, N$  trois objets de
$\bf T$. Fixons deux suites exactes de $\bf T$~:
$$0 \mapright{} A \mapright{j}
      M_1 \mapright{\pi} N \mapright{}
      0 ~, ~~~0 \mapright{} N \mapright{\iota}
      M_2 \mapright{\varpi} B \mapright{}
      0 ~.$$

Une {\it extension
panach\'ee} de $M_2$ par $M_1$ est, selon Grothendieck
([G], \S 9.3), la donn\'ee d'un diagramme commutatif de suites
exactes:
$$
\def\mapright#1{\smash{
    \mathop{\longrightarrow}\limits^{#1}}}
\def\mapdown#1{\Big\downarrow
    \rlap{$\vcenter{\hbox{$\scriptstyle#1$}}$}}
\matrix{&&&&0&&0\cr
    &&&&\mapdown{}&&\mapdown{}\cr
    0&\mapright{}&A&\mapright{j}&
      M_1&\mapright\pi&N&\mapright{}&0\cr
    &&\Big\Vert&&\mapdown{\tilde \iota}&&\mapdown\iota\cr
    0&\mapright{}&A&\mapright{\tilde j}&
      M&\mapright{\tilde \pi }&
      M_2&\mapright{}&0\cr
&&&&\mapdown{\tilde \varpi}&&\mapdown{\varpi}\cr
    &&&&B&=&\hidewidth B \hidewidth\cr
    &&&&\mapdown{}&&\mapdown{}\cr
    &&&&0&&0&&,\cr}$$
les deux 1-extensions (horizontale et verticale) o\`u s'inscrit  l'objet
$M$ v\'erifiant donc : $\iota^*(M)  \simeq M_1, \pi_*(M) \simeq  M_2$. Le choix du rel\`evement $\tilde \iota$ de $\iota$ et du prolongement $\tilde \pi$ de $\pi$ fait partie de la d\'efinition de l'extension panach\'ee, m\^eme si nous la noterons abusivement $M$. 
Lorsque $\bf T$ est tannakienne, le dual $\check M$ de
$M$ est naturellement muni d'une structure d'extension
panach\'ee $(\check M, ~^t \tilde \varpi, ~^t \tilde j)$ de la 1-extension $\check M_1$ de $\check A$ par $ \check N$ par
la 1-extension $\check M_2$ de $\check N$ par $ \check B$.

\medskip
Un morphisme $F : M \rightarrow M'$ entre deux extensions panach\'ees de
$M_2$ par $M_1$ est un  $\bf T$-morphisme induisant l'identit\'e sur $M_1$ et sur $M_2$.  En  particulier, tout endomorphisme d'une extension panach\'ee $M$ est un automorphisme, de la forme $id_M + \tilde j \circ f \circ \tilde \varpi$ pour un unique \'el\'ement $f$ de $Hom(B,A)$. Dans le m\^eme esprit,   les fl\`eches $\psi : M \rightarrow \check M$ intervenant au \S 2 dans la d\'efinition des extensions panach\'ees autoduales induisent des $\bf T$-morphismes 
{\it fix\'es}  $\Phi : M_1 \rightarrow \check M_2, ~\varepsilon \,^t \Phi : M_2 \rightarrow \check M_1$ sur les  1-extensions de d\'epart.

\medskip

Soit $E \in Ext^2(B, A)$ le produit de Yoneda de la classe de
$M_1$ dans
$Ext^1(N,A)$ par celle de $M_2$ dans $Ext^1(B,N)$.
On peut repr\'esenter
$E$ par la suite exacte:
$$ 0 \mapright{} A \mapright{j}
      M_1 \mapright{e} M_2 \mapright{\varpi}
      B \mapright{} 0 ~,$$
o\`u  $e = \iota $o$ \pi$. Une autre fa\c con
de voir
$E$ consiste \`a consid\'erer la suite exacte :
$$ Ext^1(B,A)  \mapright{j_*} Ext^1(B,M_1) \mapright{\pi_*}
      Ext^1(B,N) \mapright{\delta_{M_1}} Ext^2(B,A)$$
attach\'ee \`a l'extension $M_1$. On a alors:
$$ E ~= ~\delta_{M_1}(M_2),$$
de sorte que $E$ est nulle si et seulement
si la classe de $M_2$ appartient au noyau de $\delta_{M_1}$,
c'est-\`a-dire \`a l'image de
$\pi_*$ , c'est-\`a-dire encore si et
seulement si $M_2$ et
$M_1$ sont ``panachables". Ainsi :

\medskip

\noindent
{\bf Lemme 1} ([G], 9.3.8.c) : {\it les deux 1-extensions $M_2, M_1$ de $\bf T$ sont
panachables si et seulement si leur produit de Yoneda est nul dans $Ext^2(B,A)$.}

\medskip
Cet \'enonc\'e fournit \`a rebours un crit\`ere commode pour v\'erifier qu'une
2-extension de $B$ par $A$, donn\'ee par une suite exacte
\`a 4 termes, est \'equivalente \`a la 2-extension triviale.

\bigskip
Comme le note Grothendieck,  l'ensemble $Extpan(M_2,M_1)$ des classes
d'isomorphisme d'extensions panach\'ees de $M_2$ par $M_1$ est naturellement muni 
 d'une action du groupe
$Ext^1(B,A)$, d\'efinie  de la fa\c{c}on suivante.
\`A isomorphisme pr\`es, le translat\'e $M * U$ d'une
extension panach\'ee
$M$ de $M_2$ par $M_1$, par une 1-extension $U$ de $B$ par
$A$, est la somme de Baer  de
$M$ et de
$j_*(U)$, vues comme des 1-extensions de  $B$ par $M_1$; puisque $\pi_*(M*U) $ = $\pi_*(M) +
\pi_*$o$j_*(U) = \pi_*(M) = M_2$ dans $Ext^1(B,N)$, et que les isomorphismes associ\'es fournissent un prolongement canonique de $\pi$ \`a $M*U$, il s'agit l\`a bien
d'une extension panach\'ee de $M_2$ par $M_1$.  

\medskip

\noindent
{\bf Lemme 2} ([G], 9.3.8.b) :  {\it pour cette action (et s'il est non vide), l'ensemble $Extpan(M_2,M_1)$ des classes
d'isomorphisme d'extensions panach\'ees de $M_2$ par $M_1$ est un torseur sous le groupe 
$Ext^1(B,A)$.}

\medskip
Nous rappelons dans un appendice comment \'etablir cette propri\'et\'e, et montrons que l'autre fa\c{c}on naturelle de faire agir $Ext^1(B,A)$ sur  $Extpan(M_2,M_1)$ conduit \`a la m\^eme construction. 
En particulier, le Lemme A.1 de l'appendice 
permettrait, sous l'hypoth\`ese de rigidit\'e d\'ecrite ci-dessous,   de munir l'ensemble
$~{\cal P}(B,N,A)$ des classes d'\'equivalence
d'extensions panach\'ees d'une 1-extension (non
sp\'ecifi\'ee) de
$B$ par $N$, par une 1-extension (non sp\'ecifi\'ee) de $N$
par $A$, d'une structure de  biextension de $Ext^1(B,N)
\times Ext^1(N,A)$ par
$Ext^1(B,A)$. Nous n'en ferons  pas usage dans ce qui suit. 

\medskip
Signalons dans une direction voisine l'ensemble $~{\cal F}(B,N,A)$ introduit (dans un cadre plus g\'en\'eral) dans [RSZ], \S 2.3.1,  pour \'etudier les classes d'isomorphisme d'objets $M$ de $\bf T$ munis d'une filtration \`a 3 crans   \`a gradu\'es $A ,  N, B$ donn\'es.  Les automorphismes de $M$ n'induisent plus n\'ecessairement l'identit\'e sur les extensions interm\'ediaires $M_1, M_2$, de sorte que des extensions panach\'ees non isomorphes peuvent fournir des objets filtr\'es isomorphes.  Mais si 
$$Hom(N, A) = 0 ~, ~Hom(B, N) = 0 ~ , $$
 les 1-extensions $M_1$ et $M_2$ n'ont pas d'automorphismes non triviaux, et on peut alors identifier $~{\cal P}(B,N,A)$ et $~{\cal F}(B,N,A)$. Pour simplifier l'expos\'e, c'est sous cette   hypoth\`ese de rigidit\'e   que nous nous pla\c {c}ons maintenant.
 Voir \'egalement [BK], E.3.2. 
\bigskip

\noindent
{\bf \S 2. Extensions panach\'ees autoduales.}

\bigskip

On suppose d\'esormais que
$\bf T$ est une {\it cat\'egorie tannakienne neutre sur un corps $k$
de caract\'eristique nulle}, dont on note ${\bf 1} = \check
{\bf 1}$ l'objet neutre, $\check .$ la dualit\'e,  $^t.$
  la  transposition
sur les morphismes (et en abr\'eg\'e: $Ext^1 = Ext$). Tout objet de $\bf
T$ s'identifie canoniquement \`a son bidual, et tout morphisme \`a son
bitranspos\'ee (cf. [DM], 1.7). On se propose d'\'etudier, dans
$\bf T$, les extensions panach\'ees qui sont ``autoduales en un sens
panach\'e".

\medskip

    Plus pr\'ecis\'ement, partons
d'objets
$A$, $B$, $N$ et de 1-extensions
$M_1$,
$M_2$ comme au \S 1, pour lesquels  on suppose dans toute la
suite que
$$A \simeq \check B ~, ~N \simeq \check N~,~Hom(N,A)=  0  \qquad \qquad (0)$$
(de sorte que $Hom(B,  N) \simeq Hom(B, \check  N) $ est \'egalement nul),   et {\it fixons} deux
isomorphismes
$\phi : N \rightarrow \check N$
et $\lambda : A \rightarrow \check B$, ainsi qu'un signe $\varepsilon \in \{-1, 1\}$,  tels que

\medskip

$\bullet$  {\it le transpos\'e de $\phi$ v\'erifie
   $^t\phi = \varepsilon \phi$} ;

$\bullet$ {\it les extensions $\lambda_*M_1$ et $\phi^*\check M_2$
sont \'egales dans}
$Ext(N, \check B)$ (ou encore, par transposition,
$(^t\phi)_*M_2 = (^t\lambda)^*\check M_1$ dans $Ext(B, \check
N)$); autrement dit,
il existe un rel\`evement $\Phi : M_1
\rightarrow
\check M_2$ de $\phi$, {\it unique} d'apr\`es l'hypoth\`ese de rigidit\'e faite sur
$M_1$ et $M_2$, induisant $\lambda$ sur $A$ :

$$
\def\mapright#1{\smash{
    \mathop{\longrightarrow}\limits^{#1}}}
\def\mapdown#1{\Big\downarrow
    \rlap{$\vcenter{\hbox{$\scriptstyle#1$}}$}}
\matrix{
    0&\mapright{}&A&\mapright{j}&
      M_1&\mapright{\pi}&N&\mapright{}&0\cr
    &&\mapdown{\lambda}&&\mapdown{\Phi}&&\mapdown\phi&&&& (1)\cr
    0&\mapright{}&\check B&\mapright{^t \varpi}&
      \check M_2&\mapright{^t \iota}&
     \check N&\mapright{}&0.\cr}$$

\medskip

Supposons enfin {\it $M_2$ et $M_1$ panachables}, et soit $M$ une
extension panach\'ee
  de
$M_2$ par
$M_1$. On dira
que $M$ est {\it autoduale  relativement \`a $\Phi$} (ou,
  moins pr\'ecis\'ement, autoduale au sens panach\'e)
si $\Phi$ s'\'etend en un morphisme $\psi$ :

$$
\def\mapright#1{\smash{
    \mathop{\longrightarrow}\limits^{#1}}}
\def\mapdown#1{\Big\downarrow
    \rlap{$\vcenter{\hbox{$\scriptstyle#1$}}$}}
\matrix{
    0&\mapright{}&M_1&\mapright{}&
      M&\mapright{}&B&\mapright{}&0\cr
    &&\mapdown{\Phi}&&\mapdown{\psi}&&\mapdown{\mu = \varepsilon
^t\lambda}&&&& (2)\cr
    0&\mapright{}&\check M_2&\mapright{}&
      \check M&\mapright{}&
     \check A&\mapright{}&0\cr}$$

\noindent
  de $M$ sur $\check M$, tel que {\it le morphisme $\mu$ induit par $\psi$
sur $B$ par passage au quotient soit \'egal \`a} $
\varepsilon ^t
\lambda$  \footnote{$^{(2)}$}
{Cette deuxi\`eme
condition d\'ecoule automatiquement de la premi\`ere si
$Hom(M_1,A)$ $ = Hom(B, M_2) = 0$ (cf. Note $(3)$ plus bas), mais pas en
g\'en\'eral (penser au cas d'extensions scind\'ees). Nos conditions sont  n\'eanmoins compatibles \`a la d\'efinition donn\'ee dans [By], 2.1.4 et [M], 4.2.6, de polarisation sur un 1-motif, en ce sens qu'un 1-motif est polarisable si et seulement si l'extension panach\'ee qui lui correspond est autoduale, avec $\varepsilon = -1$. Voir la Note (4) ci-dessous pour une justification de ce signe.}. De fa\c {c}on
g\'en\'erale, on notera
$Extpanaut(M_2,M_1,\Phi)$ l'ensemble des classes d'isomorphismes
d'extensions panach\'ees
  $M$ de $M_2$ par $M_1$ qui sont autoduales relativement \`a  $\Phi$. On
verra au Lemme 4 ci-dessous que $M$ est dans ce cas automatiquement munie
d'au moins une autodualit\'e $\psi = \varepsilon ^t\psi: M
\rightarrow
\check M$ de m\^eme signe que $\phi$.

\medskip

Notant $Ext_{\pm}(B,\check B)$ le groupe des classes d'isomorphisme d'extensions
  $C$ de $B$ par $\check B$ telles que
$\check C \simeq \pm C$, on peut alors
pr\'eciser la Proposition 9.3.8.b de [G] de la fa\c{c}on
suivante.

\medskip
\noindent
{\bf Th\'eor\`eme 1} : {\it $Extpanaut(M_2,M_1, \Phi)$ est
non vide, et est de fa\c {c}on naturelle  un torseur
sous le groupe $Ext_{\varepsilon}(B,\check B)$. En particulier, il
  co\"{\i}ncide avec $Extpan(M_2,M_1)$ si  $Ext(B, \check B)
= Ext_{\varepsilon}(B,\check B)$, et est r\'eduit \`a un \'el\'ement
si $Ext(B, \check B)  = Ext_{- \varepsilon}(B,\check B)$.}

\medskip
\noindent
{\it D\'emonstration} : soit $M$ une extension panach\'ee de $M_2$
par $M_1$. L'existence
du prolongement
$\psi$ de $\Phi$ d\'ecrit par le diagramme (2)
\'equivaut, apr\`es transposition,  \`a
l'\'egalit\'e
$$\lambda_*M = \varepsilon (^t \Phi)^*
\check M ~{\rm ~dans} ~ Ext(M_2, \check B).$$
Mais les deux extensions en question sont en fait des
extensions panach\'ees de  $M_2$ par $\phi^*
\check M_2$~: c'est clair pour la seconde, et cela r\'esulte pour la premi\`ere de l'identification canonique $\lambda_* M_1 \mapright{\sim} \phi^*
\check M_2$ fournie par nos hypoth\`eses (voir (0) et (1)). Il existe
donc  un unique
\'el\'ement
$$\gamma_M = \gamma_{M,\Phi} \in Ext(B, \check B)$$
tel que $$ \lambda_*M =
\varepsilon (^t \Phi)^* \check M ~   *  \gamma_M$$
dans le $Ext(B, \check B)$-torseur
  $Extpan(M_2, \phi^* \check M_2) $.

\medskip

	Cette extension $\gamma_M$ de $B$ par $\check B$ mesure
ainsi  {\it l'obstruction \`a
l'autodualit\'e} de $M$ relative
\`a $\Phi$. Sa	 duale dans $\bf T$ v\'erifie :

\medskip

\noindent
{\bf Lemme 3} : {\it sous les hypoth\`eses du diagramme ~{\rm(1)}, on
a dans $Ext(B, \check B)$ :
$$\check {(\gamma_M) }= - \varepsilon \gamma_M.$$}
\noindent
{\it D\'emonstration du Lemme} 3 : consid\'erons la suite exacte
$$ 0 = Hom(N, \check B) ~ \mapright{(.)_*M_2} ~ Ext(B, \check B)
~ \mapright{\varpi^*} ~Ext(M_2,
\check B)
~\mapright{\iota^*} ~ Ext(N,
\check B)$$
attach\'ee \`a $M_2$. Les extensions $\lambda_*M$ et $
\varepsilon (^t \Phi)^*
\check M$ se projetant
respectivement sur
$\lambda_*M_1$ et $\varepsilon (^t \phi)^* \check M_2 = \phi^*\check M_2$, qui
sont
\'egales par hypoth\`ese dans $Ext(N, \check B)$,
 l'injectivit\'e ici suppos\'ee de $\varpi^*$ permet de d\'efinir $\gamma_M$ comme l'unique \'el\'ement  $\gamma = \gamma_M $ de $Ext(B, \check B)$ tel que
$$\lambda_*M - \varepsilon (^t
\Phi)^*
\check M = \varpi^*\gamma ~~{\rm dans}~~ Ext(M_2, \check B).
\qquad \qquad (3)$$ Plus pr\'ecis\'ement, nos extensions sont
  des extensions panach\'ees de
$M_2$ par $\phi^* \check
M_2 $
:
$$
\def\mapright#1{\smash{
    \mathop{\longrightarrow}\limits^{#1}}}
\def\mapdown#1{\Big\downarrow
    \rlap{$\vcenter{\hbox{$\scriptstyle#1$}}$}}
\matrix{&&&&0&&0\cr
    &&&&\mapdown{}&&\mapdown{}\cr
    0&\mapright{}&\check B&\mapright{^t{\hat \varpi}}&
      \{ \lambda_* M_1 \mapright{\sim}  \phi^* \check M_2\} &\mapright{}&N&\mapright{}&0\cr
    &&\Big\Vert&&\mapdown{}&&\mapdown{}\cr
    0&\mapright{}&\check B&\mapright{}&
      \lambda_*M ~{\rm resp.} ~ \varepsilon (^t \Phi)^* \check M &\mapright{}&
      M_2&\mapright{}&0\cr
&&&&\mapdown{}&&\mapdown{\varpi}\cr
    &&&&B&=&\hidewidth B \hidewidth\cr
    &&&&\mapdown{}&&\mapdown{}\cr
    &&&&0&&0&&,\cr},$$
o\`u $ ^t{\hat \varpi}$ d\'esigne le morphisme 
de
$\check B$ dans $\phi^* \check M_2$ sous-jacent au diagramme (1), et l'on a :
$$\lambda_*M = \varepsilon
(^t
\Phi)^*
\check M * \gamma ~~{\rm dans}~~
Extpan(M_2, \phi^* \check M_2).$$
D'apr\`es le lemme A.1 de l'appendice, la relation (3) s'\'ecrit
donc aussi :
$$\lambda_*M - \varepsilon
(^t
\Phi)^*
\check M = ~ ^t{\hat \varpi}_*\gamma ~~{\rm dans} ~Ext(B, \phi^* \check M_2).$$
Poussons cette derniere relation dans $Ext(B, \check M_2)$ par le morphisme
canonique $$\hat \phi: \phi^*\check M_2 \rightarrow \check
M_2
$$
(avec $\hat \phi \circ  \,^t\hat \varpi = ~^t\varpi : \check B
\rightarrow \check M_2$.)
Comme $\lambda_*M \in  Ext(B,
\lambda_* M_1)~ \mapright{\sim} ~Ext(B, \phi^* \check M_2)$ se pousse sur $\Phi_* M \in Ext(B, \check M_2)$,
$^t \Phi ^* \check M$ sur $^t \lambda^* \check M$, et
$^t{\hat
\varpi}_*\gamma$ sur
$^t \varpi_* \gamma$, on obtient:
$$\Phi_* M - \varepsilon  ^t \lambda^* \check M = ~^t\varpi_*\gamma ~{\rm
dans}~ Ext(B, \check M_2),$$
d'o\`u par dualit\'e:
$$  \lambda_*  M - \varepsilon ^t \phi^* \check M   = - \varepsilon
\varpi^*
\check \gamma ~{\rm
dans}~ Ext(M_2, \check B).$$
Cette derni\`ere relation, jointe \`a (3) et \`a
l'injectivit\'e de $\varpi ^*$, montre bien que $\check
\gamma = - \varepsilon \gamma$.

\medskip
\noindent
{\it Fin de la preuve du Th\'eor\`eme 1} : pour
  $\delta \in Ext(B, A)$, posons $\delta' = \lambda_*\delta
\in Ext(B, \check B)$. Montrons d'abord que si on
remplace
$M \in
Extpan(M_2,M_1)$ par $M *
\delta$, l'obstruction $\gamma_M$ \`a l'autodualit\'e de $M$
  devient
$$ \gamma_{M * \delta} =
\gamma_M + \delta' - \varepsilon \check \delta' \in Ext(B,
\check B).$$ Dans ce passage,  $M$ devient en effet $M
+ \varpi^*\delta$ dans $Ext(M_2,A)$, et $\lambda_*M$ donne
$\lambda_* M + \varpi^* \lambda_* \delta$ dans $Ext(M_2,
\check B)$. Comme $M$ devient $M + j_* \delta$ dans
$Ext(B,M_2)$, $\check M$ donne $\check M + {^t j}^* \check
\delta$ dans $Ext(\check M_2, \check B)$, et
$^t {\Phi}^*
\check M$ devient $ ^t {\Phi}^*
\check M + \varpi^*(\check{\lambda_* \delta})$ dans
$Ext(M_2, \check B)$ (rappelons que $\Phi {\rm o} j=$
$^t\varpi {\rm o} \lambda$, cf. (1)). Finalement,
$\varpi^* \gamma_M = \lambda_*  M - \varepsilon ^t \phi^* \check
M $ devient $\varpi^* \gamma_M + \varpi^* \lambda_* \delta -
\varepsilon \varpi^* \check{(\lambda_* \delta)}$, et
$\gamma_M$ est ainsi bien remplac\'e par $\gamma_M +  \delta'
  - \varepsilon \check{\delta'}$.

Faisons maintenant agir les \'el\'ements $\delta' = \varepsilon
\check \delta'$ de  $Ext_{\varepsilon}(B,\check B)$ sur les
objets
$M \in Extpan(M_2,M_1)$ en posant $M + \delta' = M*
\lambda_*^{-1}\delta'$. Si $M \in Extpanaut(M_2,M_1,
\Phi)$, i.e. si
$\gamma_M = 0$, il en sera de m\^eme de l'obstruction
$\gamma_{M *\delta} = \gamma_M$ \`a l'autodualit\'e de $M +
\delta'$. Inversement, tout objet $M'$ de $Extpanaut(M_2,M_1,
\Phi)$ s'\'ecrit de fa\c {c}on unique $M * \delta$
dans $Extpan(M_2,M_1)$, o\`u $\lambda_*\delta:= \delta'$
v\'erifie: $\delta' - \varepsilon \check{\delta'} = \gamma_{M'}
-
\gamma_M = 0$, de sorte qu'on a bien alors $\delta' \in
Ext_{\varepsilon}(B,
\check B)$ et $M' = M + \delta'$.

Enfin, $Extpanaut(M_2,M_1, \Phi)$ est non vide : pour tout $M \in
Extpan(M_2,M_1)$, l'obst-ruction
$\gamma_M \in Ext_{-\varepsilon}(B,\check B)$ peut d'apr\`es le Lemme 3
s'\'ecrire
$ \varepsilon \check \delta' - \delta'$, avec $\delta' = -{1\over
2}\gamma_M$ (la division par 2 a un
sens, puisque ${1\over 2} \in k \subset End(B)$); rempla\c {c}ant
$M$ par $M' = M * \lambda_*^{-1}\delta'$, on obtient $\gamma_{M'} =
0$, d'o\`u un \'el\'ement $M'$ dans $Extpantaut(M_2,M_1,\Phi)$.  

\bigskip

  On peut par ailleurs pr\'eciser de combien de fa\c{c}ons,
et avec quel signe, $\Phi$ se prolonge de fa\c {c}on panach\'ee
aux
\'el\'ements
$M$ de
$Extpanaut(M_2,M_1,\Phi)$. Posant
$Hom_{\pm }(B, \check B) = \{v \in Hom(B, \check
B), ^tv = \pm  v\}$, on obtient :

\medskip

\noindent
{\bf Lemme 4} : {\it  Soit $M \in
Extpanaut(M_2,M_1, \Phi)$. L'ensemble
$Isoaut_{\varepsilon}(M,
\Phi)$ des isomorphismes
$\psi: M
\rightarrow
\check M$ prolongeant de fa\c {c}on panach\'ee $\Phi$ avec la parit\'e
de $\phi$ (i.e. tels que $^t\psi = \varepsilon \psi$) est
non vide, et est naturellement muni d'une struture de torseur
sous le groupe
$Hom_{\varepsilon}(B, \check B)$.}

\medskip

\noindent
{\it D\'emonstration} :  montrons d'abord que pour tout
prolongement panach\'e $\psi$ de $\Phi$ \`a $M$, il existe un
\'el\'ement
$v$ de
$Hom_{-\varepsilon}(B,
\check B)$ tel que
$$\psi = \varepsilon ^t \psi + ~^t \tilde \varpi v \tilde
\varpi \in Hom(M, \check M).$$
En effet, le diagramme (2)
induit :
$$
\def\mapright#1{\smash{
    \mathop{\longrightarrow}\limits^{#1}}}
\def\mapdown#1{\Big\downarrow
    \rlap{$\vcenter{\hbox{$\scriptstyle#1$}}$}}
\matrix{
    0&\mapright{}&A&\mapright{\tilde j}&
      M&\mapright{}&M_2&\mapright{}&0\cr
    &&\mapdown{\lambda}&&\mapdown{\psi}&&\mapdown{\overline \psi}\cr
    0&\mapright{}&\check B&\mapright{^t \tilde \varpi}&
      \check M&\mapright{}&
     \check M_1&\mapright{}&0 ~,  \cr}  $$
ainsi que, par passage au quotient  :
$$
\def\mapright#1{\smash{
    \mathop{\longrightarrow}\limits^{#1}}}
\def\mapdown#1{\Big\downarrow
    \rlap{$\vcenter{\hbox{$\scriptstyle#1$}}$}}
\matrix{
    0&\mapright{}&N&\mapright{}&
      M_2&\mapright{}&B&\mapright{}&0\cr
    &&\mapdown{\phi}&&\mapdown{\overline \psi}&&\mapdown{\mu}\cr
    0&\mapright{}&\check N&\mapright{}&
      \check M_1&\mapright{ }&
     \check A&\mapright{}&0.\cr}$$
Comme l'extension $M_1$ n'a pas d'automorphisme et que $^t \mu =
\varepsilon \lambda$ par d\'efinition\footnote{$^{(3)}$}
{Sous la seule hypoth\`ese que $\psi$ prolonge $\Phi$, on obtient
   $(\lambda - \varepsilon ^t \mu)_*M_1 = (\phi - \varepsilon ^t
\phi)^*\check M_2
= 0$ dans
$Ext(N, \check B)$. Consid\'erant la suite exacte
$ Hom(M_1, \check B) \rightarrow Hom(A, \check B) \rightarrow Ext(N,
\check B)$ attach\'ee
\`a l'extension $M_1$, on en d\'eduit que la condition de quotient $ \mu =
\varepsilon\, ^t\lambda$ est automatiquement satisfaite
si $Hom(M_1,A) \simeq Hom(M_1, \check B)$ est r\'eduit \`a 0 (cf.  Note
(2)).} 
de $\psi$, la comparaison de
  ce dernier diagramme au transpos\'e de (1) donne $  \overline \psi =  \varepsilon \,^t \Phi$. Ainsi, $\varepsilon \psi^{-1} ~ ^t \psi$ est un automorphisme de
l'extension panach\'ee $M$,
donc de la forme $id_M + \tilde j u \tilde
\varpi$ pour un
\'el\'ement
$u$ de
$Hom(B, A)$. Comme $\psi \tilde j = ~^t \tilde \varpi \lambda$, on conclut en
posant
$v =
-\lambda u$, qui appartient bien \`a $Hom_{-\varepsilon}(B,
\check B)$ puisque $^t \tilde \varpi (^tv + \varepsilon v) \tilde \varpi = 0$.

Dans ces conditions, l'isomorphisme $\psi' := \psi - ^t
\tilde
\varpi ({1\over 2}v)
\tilde
\varpi = \varepsilon ^t\psi': M \rightarrow \check M$ prolonge encore
$\Phi$ de fa\c{c}on panach\'ee. Donc
$Isoaut_{\varepsilon}(M,
\Phi)$ est non vide, et on v\'erifie par le m\^eme
argument que supra que la loi $(\psi, w) \rightarrow \psi
+ ^t \tilde \varpi w \tilde \varpi$ en fait un torseur sous
$Hom_\varepsilon(B, \check B)$.

\medskip

En liaison avec la description cohomologique qu'en donne B. Kahn dans [K],
A.15, 
calculons pour terminer les groupes structuraux
$Ext_{\pm}(B, \check B)$ du th\'eor\`eme 1. \'Etant donn\'e deux objets
$A, B$ de la cat\'egorie (rigide) $\bf T$ et une extension $C$ de $B$ par
$A$, notons $F(C)$ l'extension de $\bf 1$ par $\check B \otimes A \simeq 
{\underline {Hom}}(B, A)$, image inverse de l'extension 
$0 \rightarrow
\check B \otimes A \rightarrow \check B \otimes C \rightarrow \check B
\otimes B \rightarrow 0$
sous le morphisme naturel de $\bf 1$ dans $\check B \otimes B$. Pour
tout couple $A,B$, on obtient ainsi un isomorphisme
 $$F : Ext(B, A)
~\mapright{\sim} ~ÊExt({\bf 1},
\check B\otimes A) .$$ D\'esignons encore par $t : \check B \otimes A
\rightarrow A
\otimes \check B$ la contrainte de commutativit\'e, identifi\'ee par
l'isomorphisme de bidualit\'e $({\check A})^{\check .} ~ \simeq A$ \`a
la transposition sur ${\underline {Hom}}(B, A)$, et, pour $A = \check B$
et $\varepsilon = \pm 1$, par
$\otimes ^2_\varepsilon \check B$ le noyau de $t - \varepsilon ~id$  sur
$\otimes ^2
\check B$ (autrement dit, $\otimes ^2 _{-} = \Lambda^2, \otimes ^2 _{+} =
S^2$).

\medskip
\noindent
{\bf Lemme 5} : {\it i) pour tout  $C \in Ext(B, A)$, de duale
$\check C \in Ext(\check A, \check B)$, on a
$$F(\check C) = - t_* F(C) ~Ê{\rm dans}~ Ext({\bf 1}, A \otimes
\check B)~;$$
ii) en particulier, pour tout $B \in {\bf T}$ et tout $\varepsilon = \pm
1$,
$F$ induit un isomophisme de $Ext_\varepsilon(B, \check B)$ sur $Ext({\bf
1},
\otimes_{-\varepsilon}^2 \check B)$.}

\medskip
\noindent
{\it D\'emonstration}: i) notons $i , p$ les morphismes sous-jacents \`a
l'extension $C$. On a, en identifiant $\bf 1$ \`a son image dans $\check B
\otimes  B$ (et, pour simplifier l'\'ecriture, les objets de $\bf T$ \`a des modules sur une $k$-alg\`ebre $R$) :
$$F(C) = \{f \in \underline{Hom}(B, C), ~ \exists ~\alpha:= \tilde p(f)
\in   R, ~  p\circ f = \alpha ~id_{B} \}.$$
De m\^eme, $F(\check C) = \{^tg \in  \underline{Hom}(\check A,
\check C), 
~\exists
~\beta := \tilde p'(^tg) 
\in R, ~  ^ti\circ  \,^tg~= ~^t(g\circ i) = \beta id_{\check A}\}.$
Consid\'erons alors le morphisme de $F(C)$ dans $F(\check C)$ qui envoie
$f$ sur le transpos\'e de $g := \tilde p(f) id_C - f \circ p : C
\rightarrow A$. Il induit sur le noyau $\underline{Hom}(B,A)$ de
$\tilde p$ l'application $ f \mapsto ~^t(-f) \in \underline{Hom}(\check A, \check B) =
Ker (\tilde  p')$, et par passage aux quotients, l'application identit\'e
sur $\bf 1$, puisque $(\alpha ~id_C - f \circ p) \circ i = \alpha ~id_A$.
Par cons\'equent, $ F(\check C) = (t \circ [-1])_*F(C) = t_*([-1]_*F(C)) =
-t_*F(C)$ dans $Ext({\bf 1}, A \otimes
\check B)$.

ii) $F$ d\'efinit un isomorphisme de $Ext(B, \check B)$
sur $ÊExt({\bf 1},
\otimes^2 \check B)$,  et on   d\'eduit de i) que $\check C \simeq
\varepsilon C$ si et seulement si 
$t_*(F(C)) \simeq  -
\epsilon F(C)$.

\medskip

Sous les
hypoth\`eses du diagramme (1), on peut donc \'enoncer, en paraphrase
 du th\'eor\`eme 1 : {\it

\noindent
- si
$Ext({\bf 1},  \otimes_{\varepsilon}^2 \check B) = 0$,
toute extension
panach\'ee de $M_2$ par $M_1$ est autoduale ;

\noindent
- si
$Ext({\bf 1},  \otimes_{-\varepsilon}^2 \check B) = 0$, il existe une
unique extension panach\'ee autoduale de $M_2$ par $M_1$.}

\bigskip

\noindent
{\bf \S 3 Application aux repr\'esentations unipotentes}

\medskip

Soient $M$ une extension panach\'ee dans la cat\'egorie
tannakienne neutre $\bf T$, $\omega$ un foncteur fibre sur $k$, et $G_M/k
$ le groupe alg\'ebrique \`a travers lequel le sch\'ema en groupes
${\underline {Aut}}_k^\otimes(\omega)$ agit sur le $k$-espace vectoriel
$\omega(M)$.
On reprend les hypoth\`eses du
d\'ebut du \S 2, et on suppose que
$$M \in
Extpanaut(M_2,M_1, \Phi).$$
On fixe l'un des
prolongements panach\'es et $\varepsilon$-sym\'etriques de
$\Phi$
\`a
$M$, soit $\psi$, dont l'existence est assur\'ee par le lemme 4. Alors,
$\psi$ d\'efinit sur $\omega(M)$ une forme bilin\'eaire non
d\'eg\'en\'er\'ee $\varepsilon$- sym\'etrique, qu'on notera
encore $\psi$, et $\omega(A) \i \omega(M)$ est un sous-espace totalement
isotrope pour  $\psi$, dont l'orthogonal dans $\omega(M)$
s'identifie \`a $\omega(M_1)$ (puisque $\psi$ envoie $\omega(M_1)$ (resp.
$\omega(A)$) sur l'espace $\omega(\check M_2)$ (resp. $\omega(\check B)$) des
\'equations  de
$\omega(A)$ (resp. $\omega(M_1$)) dans $\omega(M)$). Autrement dit, la
filtration
naturelle $ A \i M_1 \i M$ de l'extension panach\'ee $M$
fournit la filtration standard de l'espace $\varepsilon$-sym\'etrique
$(\omega(M),
\psi)$
$$W_{\psi} :  \omega(A) \i \omega(A)^{\perp} = \omega(M_1) \i \omega(M)$$
attach\'ee \`a son sous-espace totalement isotrope $\omega(A)$. Le
groupe
$G_M$ est donc contenu dans le sous-groupe parabolique
$P_{\omega(A)}:= Aut_{W_\psi}(\omega(M))$ du groupe orthogonal ou symplectique
$Aut_{\psi}(\omega(M))$ attach\'e \`a cette filtration. Les 
radicaux unipotents
$W_{-1}G_M \i  W_{-1}P_{\omega(A)}$ de ces groupes sont form\'es des
automorphismes induisant l'identit\'e sur le gradu\'e de
$W_\psi$.

\medskip

Le sous-groupe normal $ W_{-2}P_{\omega(A)} = \{g \in P_\omega(A), (g -
id)(\omega(M)) \i \omega(A)\}$ de $W_{-1}P_{\omega(A)}$ intersecte
$W_{-1}G_M$ suivant un sous-groupe  $W_{-2}G_M$, que nous allons
maintenant d\'eterminer, en supposant que $W_{-1}G_M/W_{-2}G_M \i
W_{-1}P_{\omega(A)} / W_{-2}P_{\omega(A)}
\simeq Hom(\omega(B), \omega(N))$ (cf. Lemme 6 ci-dessous)
est aussi gros que possible. Pour tout $\delta' \in Ext(B, \check
B)$, on pose, comme dans la preuve du theor\`eme 1~:  $M +
\delta' :=   M * \lambda_*^{-1}\delta'$, et on note $W_{-2}G_{\delta'} =
W_{-1}G_{\delta'}$ le radical unipotent du groupe alg\'ebrique
$G_{\delta'}$.

\medskip

\noindent
{\bf Th\'eor\`eme 2} : {\it Soit $M$ une  extension panach\'ee
de $M_2$ par $M_1$, munie d'un isomorphisme $\Phi: M_1 \rightarrow \check
M_2$ v\'erifiant les hypoth\`eses du diagramme  {\rm (1)} du \S 2,
et telle que
$W_{-1}G_M/W_{-2}G_M = Hom(\omega(B), \omega(N))$.

i) Si $M \in
Extpanaut(M_2,M_1, \Phi)$, alors
  $W_{-2}G_M = Hom_{-\varepsilon}(\omega(B), \omega(\check B)).$

ii) De fa\c {c}on g\'en\'erale, il existe un unique \'el\'ement $\delta_M$
de $Ext_{-\varepsilon}(B, \check B)$ tel $M' = M + \delta_M$ appartienne
\`a $Extpanaut(M_2, M_1, \Phi)$. Alors, $W_{-2}G_{M'}= 
Hom_{-\varepsilon}(\omega(B), \omega(\check B)) $,
$W_{-2}G_{\delta_M}$ est contenu dans  $Hom_{\varepsilon}(\omega(B),
\omega(\check B))
$,  et
  $W_{-2}G_M \simeq W_{-2}G_{M'} \oplus W_{-2}G_{\delta_M}
$.}

\medskip

La d\'emonstration du th\'eor\`eme 2 repose sur deux ingr\'edients.
  Le premier, qui est classique (voir par exemple [Bo], p. 16, ou [Sh], Lemme 2.10), donne
la structure du parabolique
$P_{\bf A}$
attach\'e \`a un sous-espace totalement
isotrope $\bf A$ d'un espace $\varepsilon$-sym\'etrique $({\bf M},
\psi)$.
Si $\phi$ d\'esigne la forme $\varepsilon$-sym\'etrique
non d\'eg\'en\'er\'ee
induite par $\psi$ sur le quotient
$\bf N =  A^{\perp}/A$, et $\mu$ l'isomorphisme de $\bf B = M/A^{\perp}$
sur $ \bf \check A$ induit par $\psi$, on obtient :

\medskip
\noindent
{\bf Lemme 6} : {\it avec les notations pr\'ec\'edentes,

\medskip

\centerline {$W_{-2}P_{\bf A} = \{\tilde z \in Hom({\bf B, A}), ~^t\mu
{\rm{o}} \tilde z
:= z
  \in Hom_{-\varepsilon}({\bf B, \check B})\},$}

\medskip
\noindent
et le radical unipotent 
$W_{-1} P_{\bf A}$  de $P_{\bf A}$ est isomorphe au produit
semi-direct du groupe $Hom(\bf B, N)$ par le groupe
$ Hom_{-\varepsilon}(\bf B, \check B)$, muni de la loi

\medskip

\centerline {$(z, \nu)(z', \nu') = (z + z' +
{1\over 2}( \phi(\nu, \nu') -
\phi(\nu', \nu)), \nu + \nu').$}}

\medskip

Par $\phi(\nu, \nu') :{\bf B} \rightarrow \check {\bf B}$, j'entends le
morphisme 
$ b \mapsto \phi(\nu, \nu')(b)(.) = \phi(\nu(b), \nu'(.))$, d'o\`u pour
 $\nu = \beta \otimes
n, \nu' = \beta' \otimes n'\in
\check {\bf B} \otimes {\bf N}$~: $\phi(\nu, \nu') -
\phi(\nu', \nu) = \phi(n,n')(\beta \otimes \beta' - \varepsilon
\beta' \otimes \beta) \in \otimes^2_{-\varepsilon}\check {\bf B} $.
En particulier, l'image sous
$^t\mu^{-1}$ du sous-groupe d\'eriv\'e de $W_{-1}P_{\bf A}$, qui
est engendr\'ee par les $\phi(\nu, \nu') -
\phi(\nu', \nu)$ o\`u $\nu, \nu' \in Hom(\bf B, N)$,
remplit tout $ Hom_{-\varepsilon}(\bf B, \check B)$. Ainsi,

\bigskip

\centerline {\qquad \quad $ W_{-2}P_{\bf A} = (W_{-1}P_{\bf A})^{der}
\simeq Hom_{-\varepsilon}({\bf B, \check B})$, \qquad \quad}

\centerline
{\qquad \qquad \quad $W_{-1}P_{\bf A}/ W_{-2}P_{\bf A} =
(W_{-1}P_{\bf A})^{ab}\simeq Hom({\bf B, N})$ \qquad \quad (4) .Ê}

\medskip
\noindent
{\it D\'emonstration} : si $a =
dim {\bf A}, h = dim {\bf N} $, il existe une base de
$\bf M$ dans laquelle la matrice
repr\'esentative de la forme
$\psi$ est donn\'ee par
$$\pmatrix {0 & 0 & \varepsilon {\bf I}_a \cr
0 & {\bf J}_h &0 \cr
{\bf I}_a & 0 & 0 \cr},$$
o\`u $^t{\bf J}_h = \varepsilon {\bf J}_h$ repr\'esente la forme
$\varepsilon$-sym\'etrique $\phi$ sur
  $\bf N$, et $\varepsilon{\bf I}_a$ l'isomorphisme $\mu = \varepsilon
^t\lambda$
de $\bf B$ vers $\bf \check A$. Le parabolique
$P_{\bf A}$ est alors repr\'esent\'e par le
groupe des matrices
$$ {g \in GL_{2a+h}(k), g =
\pmatrix{\alpha & ^t\xi & \zeta \cr
0 & \sigma & \nu \cr
0 & 0 & ^t\alpha^{-1} \cr},}$$
avec $\alpha \in GL_a(k) \simeq Aut(\bf A)$,
$\nu \in M_{h,a}(k) \simeq Hom({\bf B, \bf N})$,
$\sigma \in GL_{h,{\bf J}_h}(k) \simeq
Aut_{\phi}(\bf N)$,
$^t\xi = -\alpha ^t\nu {\bf J}_h \sigma$,
et enfin $\zeta = \alpha(z - {1\over 2}^t\nu{\bf J}_h\nu),$
o\`u $z$ parcourt le groupe des matrices
$$Z = \{z \in M_{a,a}(k), ^tz = -\varepsilon z\},$$
i.e. le groupe $ \{\tilde z \in Hom({\bf B, A}),
  ~^t\mu {\rm{o}} \tilde z := z
  \in Hom_{-\varepsilon}({\bf B, \check B})\}.$
Son radical unipotent est donc donn\'e par
$$W_{-1}P_{\bf A}  \simeq \{g \in GL_{2a+h}(k), g =
\pmatrix{{\bf I}_a & ^t\xi & \zeta \cr
0 & {\bf I}_h & \nu \cr
0 & 0 & {\bf I}_a \cr}\},$$
avec $\nu \in M_{h,a}(k)$,  $\xi = -\varepsilon {\bf J}_h
\nu$,
$\zeta = z - {1\over 2} ^t\nu {\bf J}_h \nu $, tandis que
$W_{-2} P_{\bf A} \simeq Z$.

\bigskip

Le  deuxi\`eme ingr\'edient de la preuve a d\'ej\`a \'et\'e invoqu\'e
dans [B2],
  Lemma 2.1.

\medskip

\noindent
{\bf Lemme 7} : {\it soient $n$ une alg\`ebre
de Lie nilpotente, d'alg\`ebre d\'eriv\'ee $Dn$, et $g$ une sous-alg\`ebre
de Lie de $n$. Alors, $g = n$ si (et seulement si) $g/g \cap Dn = n/Dn$.}

\medskip

\noindent
{\it D\'emonstration} : on v\'erifie par r\'ecurrence que
la s\'erie centrale descendante
$$ C^1g = g \supset C^2g = Dg \supset ... \supset
C^i g = [C^{i-1}g,g] \supset ...$$ de $g$
v\'erifie: $C^i g + C^{i+1}n = C^i n$, de sorte que $C^i g
\rightarrow C^i n/C^{i+1} n$ est surjective pour tout
$i$. Mais
$C^i n = 0$ pour $i >> 0$, donc $C^i g = C^i n$ pour tout $i$.

\bigskip

\noindent
{\it D\'emonstration du Th\'eor\`eme 2} :

i) soit $n$ l'alg\`ebre de Lie du
radical unipotent $W_{-1}P_{\omega(A)}$
du groupe parabolique $P_{\omega(A)} = Aut_{W_\psi}(\omega(M))$, et
  $g$ celle de
$W_{-1}G$. D'apr\`es l'identit\'e (4),  $Dn = Lie W_{-2}P_{\omega(A)}$, et
$n/Dn \simeq Lie Hom(\omega(B), \omega(N))$. Comme
  $W_{-2}G_M = G_M \cap W_{-2}P_{\omega(A)}$,
l'hypoth\`ese de l'\'enonc\'e
  $W_{-1}G_M/W_{-2}G_M = Hom(\omega(B), \omega(N))$
revient donc \`a dire que $g/ g \cap Dn = n/Dn$.
Le lemme 7 entra\^ {\i}ne alors que $W_{-1}G_M = W_{-1}P_{\omega(A)}$,
et qu'en particulier $W_{-2}G_M =  W_{-2}P_{\omega(A)} =
Hom_{-\varepsilon}(\omega(B), \omega(\check B)).$

ii) on a vu au th\'eor\`eme 1 que $\delta_M := \delta' =  -{1\over
2}\gamma_M$ v\'erifie la condition demand\'ee. Si $M$ s'\'ecrit aussi 
$M''
-
\delta''$ avec $\delta'' \in Ext_{-\varepsilon}(B,\check B)$ et $M''$
autoduale, la structure de torseur de $Extpanaut(M_2, M_1, \Phi)$ fournit
un \'el\'ement
$\theta
\in Ext_{\varepsilon}(B,\check B)$ tel que $M'' = M' + \theta$, d'o\`u
$\theta = \delta_M - \delta''$ dans $Ext(B, \check B) = Ext_+(B, \check B)
\oplus Ext_-(B, \check B)$, et $\delta_M = \delta''$. 
Dans ces
conditions, $W_{-2}G_{\delta_M} \i Hom_{\varepsilon}(\omega(B),
\omega(\check B))$. Les images dans $Hom(\omega(B),
\omega(\check B))$ de  $W_{-2}G_{M'}$ et de $W_{-2}G_{\delta_M}$
sont donc lin\'eairement ind\'ependantes, et $W_{-2}G_M$ s'envoie
surjectivement sur chacune d'elles; en effet, toute extension
panach\'ee
$M$  de $M_2$ par $M_1$ (resp. toute extension $\delta$ de $B$
par $\check B$) d\'efinit un homomorphisme
$\xi_M$ (resp. $\xi_\delta$) du
sous-sch\'ema en groupes de 
${\underline {Aut}}_k^\otimes(\omega)$ d\'ecoup\'e par $M_1 \simeq \check
M_2$, 
\`a valeurs dans
$Hom(\omega(B),
\omega(\check B))$, dont l'image co\" {\i}ncide
pr\'ecis\'ement avec $W_{-2}G_M$ (resp.
$W_{-2}G_\delta$), et qui v\'erifie:
$\xi_{M +
\delta} =
\xi_M +
\xi_{\delta}$. Enfin,
$W_{-2}G_{M}$ contient
$(W_{-1}G_M)^{der} \simeq (W_{-1}G_{M'})^{der}$, qui coincide avec
$W_{-2}G_{M'}$ d'apr\`es la premi\`ere partie de la preuve. Donc
$W_{-2}G_{M}$ remplit tout  $W_{-2}G_{M'} \oplus
W_{-2}G_{\delta_M}$.

\medskip
En combinant le th\'eor\`eme 2 \`a la remarque donn\'ee \`a la fin du \S 2, on obtient en particulier :

\medskip

\noindent
{\bf Corollaire} : {\it sous les hypoth\`eses du th\'eor\`eme 2,

\noindent
i) si
$Ext({\bf 1},  \otimes_{\varepsilon}^2 \check B) = 0$,
alors $M$ appartient \`a
$Extpanaut(M_2,M_1, \Phi)$,  et par cons\'equent, 
$W_{-2}G_M$ $ = Hom_{-\varepsilon}(\omega(B), \omega(\check B))$ ;

\noindent
ii) si
$Ext({\bf 1},  \otimes_{-\varepsilon}^2 \check B) = 0$,
alors il existe un unique \'el\'ement $\delta_M$ de $Ext(B, A)$
tel que $M' = M*\delta_M \in Extpanaut(M_2,M_1, \Phi)$, et 
$W_{-2}G_{M} \simeq Hom_{-\varepsilon}(\omega(B), \omega(\check B))
\oplus W_{-2}G_{\delta_M}$.}

\bigskip
\noindent
{\it Applications} : 

\medskip
\noindent
1) Supposons que $B$ soit un objet  de  $\bf T$  inversible (par exemple,
que 
$A = B = \bf 1$) et que l'extension $M_1 \simeq \phi^*(\check M_2)$ n'admette de section au dessus d'aucun 
sous-objet non nul de
$N$. L'hypoth\`ese 
$W_{-1}G_M/W_{-2}G_M  \simeq Hom(\omega(B),
\omega(N)) $ du th\'eor\`eme 2 est alors satisfaite, tandis que $B$
v\'erifie trivialement :
$Ext({\bf 1}, 
\Lambda^2 \check B) = 0$. Les conclusions du corollaire au th\'eor\`eme 2
recouvrent dans ce cas
  l'ensemble des r\'esultats de [B2], \S 3, sur les \'equations
diff\'erentielles autoduales et leurs groupes de Galois. Plus
pr\'ecis\'ement:  

\medskip

si $\varepsilon
= -1$, l'extension panach\'ee $M$ est automatiquement autoduale, et
$W_{-2}G_M = Hom_+(k,k) = k$,
alors que 

si  $\varepsilon =
1$, 
$W_{-2}G_{M'} = Hom_-(k,k) = 0$, d'o\`u $W_{-2}G_M = W_{-2}G_{\delta_M}$,
qui vaut
$0$ ou
$k$ suivant que
$M$ est ou non autoduale. 

\noindent
Il serait int\'eresssant de confronter le th\'eor\`eme 2 \`a la description th\'eorique du radical unipotent des groupes de Galois diff\'erentiels obtenue par C. Hardouin [H] pour un produit de trois op\'erateurs semi-simples arbitraires. 

\medskip
\noindent
2)
En rempla\c {c}ant
$\check B = \bf 1$ et $\check N$ par des tordues \`a la Tate  (et 
en prenant garde
au signe\footnote{$^{(4)}$}
{Soient $M$ la r\'ealisation de Betti sur $\bf Q$ d'un 1-motif sur $\bf
C$,
$\check M = Hom(M, {\bf Q}(1))$ celle de son dual de Cartier, et $< , >_M:
M\otimes \check M \rightarrow {\bf Q}(1)$ l'accouplement canonique.
L'isomorphisme de bidualit\'e
$i: M
\rightarrow ({\check M})^{\check .}$ est alors donn\'e par $<x,\xi>_M
= {\bf -}<\xi,i(x)>_{\check M}$. Le transpos\'e $\check \psi$ d'un
morphisme
$\psi : M
\rightarrow \check M$ v\'erifie donc, en terme de la transposition usuelle
sur les espaces vectoriels: $\omega(\check \psi) = -^t(\omega(\psi))$.}
  qui appara\^{\i}t dans l'expression de la bidualit\'e {\it via} les
biextensions de Poincar\'e, 
 cf.  [G], 10.2.8, et [D], 10.2.4),
  on retrouve \'egalement
  les r\'esultats de K.  Ribet [R] sur les d\'eg\'en\'erescences  des repr\'esentations galoisiennes  attach\'ees aux
1-motifs de rangs  torique et constant \'egaux
\`a $\rho = 1$, et ceux de [B1] sur le radical unipotent de leurs groupes
de Mumford-Tate. L'influence de la relation $\Lambda^2 {\bf Z} = 0$ sur
ces
\'enonc\'es avait d\'ej\`a
\'et\'e remarqu\'ee par L. Breen (voir, plus g\'en\'eralement, [Br]). 

\medskip

Lorsque $\rho > 1$, l'analogue pour les 1-motifs des
hypoth\`eses (i) ou (ii) du corollaire n'est jamais satisfait. Mais il existe  des 1-motifs autoduaux de rang constant arbitraire, auquel le th\'eor\`eme 2 s'appliquera. Ainsi, tout 1-motif
polaris\'e au sens de la Note (2) (resp. ``antipolaris\'e"), construit sur
$\rho$ points $End(\cal N)$-lin\'eairement ind\'ependants d'une
vari\'et\'e ab\'elienne complexe $\cal N$,  admet $Sym^2{\bf Q}^\rho$
(resp. $\Lambda^2 {\bf Q}^\rho$) en cran $W_{-2}$ du radical
unipotent de son groupe de Mumford-Tate. Cette approche permet de retrouver, dans le cas autodual, certains des r\'esultats obtenus par C. Bertolin
[Be]  dans le cas g\'en\'eral. 

\bigskip
\centerline {\bf Appendice}

\bigskip
Nous reprenons le cadre g\'en\'eral du \S 1 relatif \`a la cat\'egorie ab\'elienne $\bf T$, dont on fixe les objets $A, B, N$ et les 1-extensions $M_1, M_2$.   Notons ${\bf Ext}(Q,P)$ la cat\'egorie des 1-extensions dans $\bf T$ de $Q$ par $P$, et ${\bf Extpan}(M_2, M_1)$ celle des extensions panach\'ees de $M_2 \in {\bf Ext}(B,N)$   par $M_1 \in {\bf Ext}(N,A)$. Comme annonc\'e dans le texte, l'action de  $Ext^1(B,A)$ sur $Extpan(M_2,M_1)$ peut se d\'ecrire des deux fa\c{c}ons \'equivalentes suivantes.

\medskip

\noindent
{\bf Lemme A.1} : {\it soit $M$ une extension panach\'ee
de $M_2$ par $M_1$ . Pour toute extension $U$ de $B$ par
$A$ , les sommes  de Baer $  M + j_*U$ dans ${\bf Ext}(B,M_1)$ ~et
~$  M + \varpi^*U$ dans ${\bf Ext}(M_2,A)$ d\'efinissent des extensions panach\'ees $M_U$ et $M^U$ de $M_2$ par
$M_1$, reli\'ees par un isomorphisme canonique dans ${\bf Extpan}(M_2,M_1)$ .}

\medskip

\noindent
{\it D\'emonstration} : les morphismes $j$ et $\varpi$ sont ceux  
du diagramme du \S 1 relatif  \`a l'extension panach\'ee
$M$, dont nous reprenons \'egalement les notations $\tilde j = \tilde \iota $o$ j$ et
$\tilde \varpi = \varpi $o$ \tilde \pi $. Soit par ailleurs 
$$0 \mapright{} A \mapright{i}
      U \mapright{p} B \mapright{}
      0 ~$$
      la suite exacte correspondant \`a $U$.  Alors (et en supposant pour faciliter la lecture que les objets de $\bf T$ sont des modules sur un anneau), $M_U := M + j_*U  ~=  M \times_B^{M_1}(U\times^A M_1)$ est repr\'esent\'e par
      $$
{\{ (m,u,m_1)\in M \times U \times M_1 ~,~\tilde \varpi(m)=p(u)\}
\over
\{(-\tilde \iota(m_1), -i(a), j(a)+m_1) ~, ~m_1
\in M_1 , a \in A \}} \simeq
~{\{ (m,u)\in M \times U ~, ~\tilde \varpi(m) =
p(u)\} \over \{(\tilde j(a), -i(a)) ~, ~a \in A \}},$$
tandis que $M^U := M + \varpi^*U = M \times_{M_2}^A(U\times_B M_2) $ est repr\'esent\'e par 
$$
{\{(m,u,m_2)\in M \times U \times M_2,\tilde \pi(m) = m_2,p(u)=\varpi(m_2)
\} \over \{(-\tilde j(a), i(a) , 0) ~, ~a
\in A \}} \simeq
{\{ (m,u)\in M \times U,p(u)=\tilde \varpi(m)\} \over
\{(-\tilde j(a), i(a)) ~, ~a \in A \}}.$$
Notons $F$ le $\bf T$-isomorphisme canonique de $M_U$ vers $M^U$  auquel l'\'egalit\'e des termes de droite conduit. 

\medskip
La structure de 1-extension de $M_U \in {\bf Ext}(B,M_1)$ est donn\'ee  par 
$$M_1 \rightarrow  M_U : m_1 \mapsto \overline{(\tilde \iota(m_1), 0)} ~;~ M_U  \rightarrow B : \overline{(m,u)} \mapsto \tilde \varpi(m)=p(u).$$
En composant le premier morphisme avec $j$, on voit que $M_U$ est naturellement muni d'une structure de 1-extension de $M_2$ par $A$,  soit $M_U \in {\bf Ext}(M_2, A)$, donn\'ee  par :
 $$A \rightarrow M_U : a \mapsto \overline{(\tilde j(a), 0)} = \overline{(0, i(a))} ~;~ M_U \rightarrow M_2 :  \overline{(m,u)} \rightarrow \tilde \pi(m).$$ 
 Ces morphismes permettent d'inscrire $M_U$ dans un diagramme de type (1), et en font donc une extension panach\'ee de $M_2$ par $M_1$ (dont la notation $M * U$ du texte repr\'esente la classe d'isomorphisme). Un calcul similaire sur $M^U$ entra\^{\i}ne alors que  $F : M_U \rightarrow M^U$ induit l'identit\'e  sur $M_1$ et sur $M_2$; c'est donc bien un isomorphisme dans ${\bf Extpan}(M_2, M_1)$. 

\bigskip

      $$
\def\mapright#1{\smash{
    \mathop{\longrightarrow}\limits^{#1}}}
\def\mapdown#1{\Big\downarrow
    \rlap{$\vcenter{\hbox{$\scriptstyle#1$}}$}}
\matrix{     &Hom(B,A)&         &     && \cr
 (5)  \qquad \qquad \qquad \qquad   &\mapdown{\varpi^*}&& &&  \cr
   &Hom(M_2,A)&         &     && \cr
  &\mapdown{\iota^*}&& &&   \cr
  &Hom(N,A)&\mapright{}&Hom(N,M_1)&&\cr
  &\mapdown{(.)_*M_2}&&\mapdown{}\cr
   Hom(B,A) ~\mapright{j_*} ~Hom(B,M_1)~\mapright{\pi_*}~Hom(B,N)  \quad  \mapright{(.)^*M_1}&
      Ext^1(B,A)&\mapright{j_*}&Ext^1(B,M_1)&&\cr
   \qquad \qquad \qquad \qquad  \qquad \qquad \qquad     \mapdown{ } &\mapdown{\varpi^*}&&\mapdown{    }\cr
 \qquad \qquad \qquad \qquad  \qquad \qquad \qquad \quad  Hom(M_2,N) \quad \mapright{ }&
      Ext^1(M_2,A)&\mapright{ }&
      Ext^1(M_2, M_1)  \hidewidth   \hidewidth\cr }$$

\medskip
En   d\'epit du  lemme A.1, la construction asym\'etrique de la classe d'isomorphisme $M*U \in Extpan(M_2, M_1) $ donn\'ee par $M_U$ (resp. par $M^U$) peut faire penser que l'action de $Ext^1(B,A)$ sur $Extpan(M_2, M_1)$ passe au quotient par l'image de $Hom(B, N)$ (resp. $Hom(N, A)$) dans $Ext^1(B,A)$ : voir la quatri\`eme ligne (resp. colonne)  du diagramme (5)   ci-dessus. Le lemme suivant (ou son analogue pour $M^U$) montre qu'il n'en est rien.     
   
\medskip
\noindent
{\bf Lemme A.2} : {\it soit $f \in Hom(B,N)$, et soit $U = f^*M_1 \in {\bf Ext}(B,A)$, de sorte que $j_*U \in {\bf Ext}(B, M_1)$  est munie d'une section  canonique $s_f$. Soit de plus $F : M  \mapright{\sim} M_U := M + j_*U$ le  ${\bf Ext}(B, M_1)$-isomorphisme attach\'e \`a $s_f$.  Alors, 

i) $F$ induit sur $M_2$ l'automorphisme $id_{M_2} - \iota \circ  f \circ \varpi$;

ii) $M$ et $M_U$ sont isomorphes dans ${\bf Extpan}(M_2, M_1)$ si et seulement si $f$ se rel\`eve en un morphisme de $B$ vers $M_1$, auquel cas la classe de $U$ dans $Ext^1(B,A)$ s'annule.}

\medskip
\noindent
{\it D\'emonstration} : i) avec les  conventions de la preuve pr\'ec\'edente, l'extension  triviale $j_*f^*M_1 = M_1 \times^A (M_1 \times_{N, f} B)$ est repr\'esent\'ee par 
$${\{ (m'_1,m_1,b)\in M_1 \times M_1 \times B ~,  \pi(m_1) = f(b)\}
\over
\{( j(a), -j(a), 0) ~, ~a \in A \}} ,$$ qui admet pour section $s_f : B \rightarrow j_*f^*M_1 : b \mapsto  \overline{(-\mu_1, \mu_1, b)}$, o\`u $\mu_1 \in M_1$ d\'esigne un \'el\'ement quelconque de la fibre de $f(b)$.
De m\^eme, $M_U = M + j_* f^*M_1 \in {\bf Ext}(B,  M_1)$ est repr\'esent\'e par 
$$
{\{ (m,m_1,b)\in M \times M_1 \times B ~,~\tilde \varpi(m)= b, \pi(m_1) = f(b)\}
\over
\{(\tilde j(a), -j(a), 0) ~, ~a \in A \}} \simeq
~{\{ (m,m_1)\in M \times M_1 ~, ~ \pi(m_1) = f(\tilde \varpi(m))\}
 \over \{(\tilde j(a), -j(a)) ~, ~a \in A \}},$$
 et l'isomorphisme $F : M \rightarrow M_U$ associ\'e \`a $s_f$ par
 $$m \mapsto F(m) = \overline{(m- \tilde  \iota( \mu_1), \mu_1)}, ~{\rm {pour ~tout}}~  \mu_1  ~{\rm tel ~que~}  \pi(\mu_1) =  f(\tilde \varpi(m)).$$
 On v\'erifie que $F$  est bien un morphisme de ${\bf Ext}(B, M_1)$. Par ailleurs, la structure d'extension de $M_2$ par $A$ que porte $M_U$ est  maintenant donn\'ee par
 $$A \rightarrow M_U : a \mapsto   \overline{( \tilde j(a),0)} = \overline{( 0, j(a))} ~;~M_U \rightarrow M_2 : \overline{(m,m_1)} \rightarrow \tilde \pi(m).$$

Par passage au quotient par $A$, $F$ induit  un ${\bf Ext}(B,N)$-endomorphisme $\tilde  F$ de $M_2$, qu'on peut expliciter comme suit. Soit $m_2$ un \'el\'ement de $M_2$, d'image $b = \varpi(m_2)$ dans $B$, et soient $\mu$ (resp. $\mu_1$) un \'el\'ement de la fibre de $M$ (resp. $M_1$) au-dessus de $b$ (resp. $f(b)$).  Alors,  $\pi(\mu_1) =  f(b) = f(\tilde \varpi(\mu))$, donc  $F (\mu) =  \overline{(\mu- \tilde  \iota( \mu_1), \mu_1)}$, et
$$ \tilde F(m_2) = \tilde \pi(\mu- \tilde \iota( \mu_1)) = m_2 - \iota(f(b)) = m_2 - \iota f \varpi (m_2).$$
Ainsi, $\tilde F  = id_{M_2} - \iota \circ  f \circ  \varpi$.

\medskip
ii) les extensions panach\'ees $M$ et $M_U$ sont isomorphes si et seulement s'il existe un $\bf T$-morphisme $F' : M \rightarrow M_U$ induisant l'identit\'e sur $M_1$ et sur $M_2$. En particulier, $F^{-1}Ê\circ F'$ est alors un   ${\bf Ext}(B, M_1)$-automorphisme de $M$, et il existe un \'el\'ement $g$ de $Hom(B, M_1)$ tel que    $F^{-1} \circ F' =  id_M +  \tilde \varpi \circ g \circ \tilde \iota$. Dans ces conditions, $F'$ r\'epond \`a la question si et seulement si le morphisme $\tilde F'$ qu'il induit sur $M_2$ est l'identit\'e.  Comme  $ \tilde \varpi \circ g \circ \tilde \iota$ induit $\varpi \circ \pi_*(g) \circ \iota$ sur $M_2$, de sorte que   $\tilde F' = id_{M_2} +   \varpi \circ (\pi_*(g) -f)  \circ \iota$,  cela revient \`a demander que $f$ appartienne \`a l'image de $\pi_*$ dans $Hom(B, N)$. La derni\`ere assertion d\'ecoule du diagramme (5). 

\medskip
On d\'eduit ais\'ement  de ce lemme que l'action de $Ext^1(B,A)$ sur $Extpan(M_2,M_1)$ est bien

\noindent
-~  libre :   soit $U \in {\bf Ext}(B,A)$ tel que les extensions panach\'ees $M$ et $M_U$ soient isomorphes.  En particulier, leurs classes dans $Ext(B, M_1)$ co\"{\i}ncident, donc il existe $f \in Hom(B,N)$ tel que $U  = f^*(M_1)$ dans $Ext^1(B,A)$,  et on conclut par le Lemme A.2.ii ;

\noindent
-~et transitive : soient $M$ et $M'$ deux  extensions panach\'ees. Leurs images dans $Ext^1(B,N)$ co\"{\i}ncident, donc il existe une 1-extension $U'$ de $B$ par $A$ telle que    $M'$ et $M_{U'}  = M + j_*U'$  soient li\'ees par un ${\bf Ext}(B, M_1)$-isomorphisme $\Phi$. Le ${\bf Ext}(B, N)$-automorphisme $\tilde \Phi$ de $M_2$ induit par $\Phi$ est de la forme $id_{M_2}Ê+ \iota \circ  \phi \circ  \varpi$, o\`u $\phi \in Hom(B, N)$, et le Lemme A.2.i entra\^{\i}ne que $M'$ est isomorphe dans ${\bf Extpan}(M_2, M_1)$ \`a $M_U$, o\`u $U = U' + \phi^*(M_1) \in {\bf Ext}(B,A)$.

 \bigskip

{\centerline {\bf R\'ef\'erences}}

\bigskip

\noindent
[BK] L. Barbieri-Viale, B. Kahn: {\it On the derived category of $1$-motives}; arXiv: math.AG 

1009.1900.

\noindent
[Be] C. Bertolin: {\it Le radical unipotent du groupe de Galois d'un $1$-motif};
Math. Ann.  {\bf 327},  

2003, 585-607. 

\noindent
[B1] D. Bertrand: {\it Relative splittings of one-motives}; Contemp. Maths,
{\bf
210}, 1998, 3--17

\noindent
[B2]  D. Bertrand: {\it Unipotent radicals of differential Galois groups};
Math. Ann. {\bf 321}, 2001,

 645-666.

\noindent
[B3] D. Bertrand: {\it   Extensions panach\'ees et dualit\'e};
 Pr\'epublications de l'Institut de Math\'e-
 
 matiques de Jussieu, No {\bf 287}, Avril 2001 (non publi\'e).

\noindent
[Bo] A. Borel: {\it Linear algebraic groups}; PSPM AMS,
vol. {\bf 9}, 1966, 3-19.

\noindent
[Br]  L. Breen: {\it Biextensions altern\'ees}; Compo. Math., {\bf 63}, 1987,
99--122.

\noindent
[By] J-L. Brylinski: {\it $1$-motifs et formes automorphes}; Publ. Math. Univ. 
Paris 7, {\bf 15}, 1983,

43-106.

\noindent
[D] P. Deligne: {\it Th\'eorie de Hodge III}; Publ. Math. IHES, {\bf 44}, 1975
5-77.

\noindent
[DM] P. Deligne, J. Milne: {\it Tannakian categories}; Springer LN {\bf 900},
1982, 101--228.

\noindent
[G] A. Grothendieck: {\it Mod\`eles de N\'eron et monodromie}; SGA VII.1, no
9, Springer LN

{\bf 288}, 1968.

\noindent
[H]Ê C. Hardouin : {\it Calcul du groupe de Galois diff\'erentiel du produit de trois op\'erateurs

 compl\`etement r\'eductibles}; CRAS Paris   {\bf 341},  2005, 349-352.

\noindent
[K] B. Kahn: {\it Repr\'esentations orthogonales et symplectiques sur un corps
de caract\'eristi-

que diff\'erente de} 2; Comm. in Algebra {\bf 31}, 2003, 133--196. 

\noindent
[M] J. Milne: {\it Canonical models of (mixed) Shimura varieties and
automorphic vector bun-

dles}; ``Automorphic forms, Shimura
varieties and $L$-functions", I, 283-411 (1990).

\noindent
[RSZ] J-P. Ramis, J. Sauloy, C. Zhang : {\it Local analytic classification of $q$-difference equa-

tions}; arXiv: math.AQ 0903.0853.

\noindent
[R] K. Ribet: {\it Cohomological realization of a family of one-motives}; J.
Number Th., {\bf 25},

1987, 152--161.

\noindent
[Sh]  G. Shimura: {\it Euler products and Eisenstein series}; CBMS AMS, vol. {\bf 93}, 1997.

\end